\documentclass{article}
\usepackage{amsmath}
\usepackage{amsfonts}
\usepackage{graphicx}
\usepackage{latexsym}
\newcommand{\qed}{\hspace*{\fill} $\Box$}
\def\endrem{}
\def\colon{{:}\;}

\newcommand{\be}[1]{\begin{equation} \label{#1}}

\newcommand{\norm}[1]{\| {#1}\| }

\def\charfct{1\hspace{-2pt} {\rm l}}

\def\schnitt{\cap}

\def\verein{\cup}
\def\Verein{\bigcup}

\def\a{\alpha}
\def\b{\beta}
\def\g{\gamma}
\def\d{\delta}
\def\e{\varepsilon}

\def\cC{{\cal C}}
\def\cT{{\cal T}}

\def\R{{\mathbb R}}

\def\N{{\mathbb N}}
\def\Z{{\mathbb Z}}

\def\P{{\mathbb P}}
\def\E{{\mathbb E}}

\newtheorem{Theorem}{Theorem}
\newtheorem{Lemma}[Theorem]{Lemma}

\begin{document}

\title{Quantitative estimates of discrete harmonic measures}
\author{E. Bolthausen and K. M\"{u}nch-Berndl\footnote{Supported by
Swiss NF grant 20-55648.98}\\ {\normalsize Institut f\"{u}r Angewandte
Mathematik der Universit\"{a}t Z\"{u}rich}\\ {\normalsize Winterthurer
Str.\ 190, 8057 Z\"{u}rich, Switzerland}}
\date{May 2, 2000}
\maketitle
\begin{abstract}
A theorem of Bourgain states that the harmonic measure for a
domain in
$\R^d$ is supported on a set of Hausdorff dimension strictly less than
$d$
\cite{Bourgain}. 
We apply Bourgain's method to the discrete case, i.e.,
to the distribution of the first entrance point of a random walk into a
subset of $\Z ^d$, $d\geq 2$. 
By refining the argument, we prove that for
all $\b>0$ there exists $\rho (d,\b )<d$ and $N(d,\b)$, such that  for
any $n>N(d,\b )$, any $x \in
\Z^d$, and any $A\subset \{ 1,\dots , n\}^d$ 
$$ 
| \{ y\in\Z^d\colon \nu_{A,x} (y) \geq n^{-\b} \}| \ \leq \ 
n^{\rho(d,\b )} ,
$$
where $\nu_{A,x} (y)$ denotes the probability that $y$ is the first
entrance point of the simple random walk starting at $x$ into $A$.
Furthermore, $\rho$ must converge to $d$ as $\b \to \infty$.
\end{abstract}

\section{Introduction} 

Let $(S_n)_{n\in \N}$ be a simple random walk in $\Z^d$ starting at
$x\in
\Z^d$, i.e., $S_0=x$ and 
$$ 
\P^x(S_{n+1} - S_n = e) = \frac 1{2d}, \quad \norm{e}=1, \quad
n\in\N.
$$ 
($\norm{ \ . \ }$ denotes the Euclidian distance, i.e.,
$\norm{x} =
\sqrt{x_1^2+\dots +x_d^2}$.) 
For
$A\subset
\Z^d$, $A\neq \emptyset$, we denote by
$\tau_A$ the time of the first entrance of $S$ to $A$:
$$ 
\tau_A=\inf \{n\geq 0\colon S_n\in A\} .
$$
The harmonic measure for $A$ of a set $B\subset \Z^d$ evaluated at
$x\in\Z^d$ is defined as
$$ 
\omega (A,B,x) = \P^x (\tau_A<\infty ,\ S_{\tau_A} \in B).
$$
Clearly, for $x\in A$, $\omega (A,B,x) = \charfct_B(x)$.
For fixed $A\subset \Z^d$ and $x\in \Z^d$, $\omega (A,\ .\ ,x)$ is a
measure on
$\Z^d$ with total mass
$\omega(A,\Z^d,x)= \omega(A,A,x) =
\P^x (\tau_A<\infty )\in (0, 1]$.  
We denote by
$\nu_{A,x}(y)=\omega (A,\{ y\}, x)$ its density. 
For
$x\in A^c =
\Z^d\setminus A$, $\omega (A,B,\ .\ )$ is a harmonic function,
$$
\Delta
\omega (A,B,x)= \frac 1{2d} \sum_{\norm{e}=1}\omega (A,B,x+e) - \omega
(A,B,x) = 0.
$$

We shall prove the following theorem:

\medskip
\noindent
{\bf Theorem.}
{\it (A) For all $\b>0$ there exists $\rho (d,\b )<d$ and $N(d,\b)$, such
that  for any $n>N(d,\b )$, any $x \in
\Z^d$, and any $A\subset Q^d(n) = \{ 1,\dots , n\}^d$ 
$$ 
| \{ y\in\Z^d\colon \nu_{A,x} (y) \geq n^{-\b} \}| \ \leq \ 
n^{\rho(d,\b )} .
$$ 
(B) For all $\rho <d$ there 
exist $\beta <\infty$ and sequences $n_K \to
\infty$,
$x_K\in \Z^d$, and $A_K\subset Q^d(n_K)$ such that for all $K$
$$ 
| \{ y\in\Z^d\colon \nu_{A_K,x_K} (y) \geq n_K^{-\b} \}| \ > \ 
n_K^{\rho} .
$$  
}

\medskip

(For $A\subset \Z^d$, $|A|$ denotes the number of points of $A$.)

\medskip
\noindent
{\bf Remarks.} (1) If $x \in A$, the statement of Theorem (A) is trivial.
Therefore we only consider $x \in A^c$.
The proof of Theorem (A) is to a large extent an adaptation 
of Bourgain's proof \cite{Bourgain} that the harmonic measure for a
domain in
$\R^d$ is supported on a set of Hausdorff dimension strictly less than
$d$ to the discrete case, 
and the proof of Theorem (B) is inspired by Jones and
Makarov \cite{JM} who also treat continuous harmonic measure. 
\endrem

(2) The analogous theorems hold for harmonic measure conditioned on the
event that $A$ is reached, and also for harmonic measure from infinity:
Let 
$$
\bar\nu_{A,x} (y) = \P^x (S_{\tau_A} =y|\; \tau_A<\infty ) ,
$$
and 
$$
\bar\nu_{A,\infty} (y) = \lim_{\| x\| \to \infty}\bar\nu_{A,x}
(y).
$$
(See for example \cite{Lawlerbook}, Chapter 2.1 for the existence of
$\bar\nu_{A,\infty}$.) 
Then we have

{\it (A') For all $\b>0$ there exists $\rho (d,\b )<d$ and $N(d,\b)$, such
that  for any $n>N(d,\b )$, any $x \in
\Z^d$, and any $A\subset Q^d(n) = \{ 1,\dots , n\}^d$ 
$$
| \{ y\in\Z^d\colon \bar\nu_{A,x} (y) \geq n^{-\b} \}| \ \leq \ 
n^{\rho(d,\b )} ,
$$ 
}
and

{\it (A'') For all $\b>0$ there exists $\rho (d,\b )<d$ and $N(d,\b)$, such
that  for any $n>N(d,\b )$ and any $A\subset Q^d(n) = \{ 1,\dots ,
n\}^d$ 
$$ 
| \{ y\in\Z^d\colon \bar\nu_{A,\infty } (y) \geq n^{-\b} \}| \ \leq \ 
n^{\rho(d,\b )} .
$$ 
}
For (A'), note first that for $d=2$, $\P^x (\tau_A<\infty ) =1$ for all
$x$ and $A$ by recurrence and therefore $\bar\nu_{A,x} = \nu_{A,x}$. 
For $d\geq 3$, we have a lower bound on the hitting probability $\P^x
(\tau_A<\infty ) $ for $x$ in a neighborhood of $Q^d(n)$, 
\begin{equation} 
\P^x
(\tau_A<\infty ) \geq \P^x
(\tau_{\{ z\} }<\infty ) = \frac {G(x-z)}{G(0)} \geq \frac {c_2}{G(0)} \|
x-z\| ^{2-d} \geq c(a,d) n^{2-d}\label{eq1}
\end{equation}
for all $z\in A$ and $x\in U^d(an) = \{ -an, \dots, (a+1)n\} ^d$, where
$G$ is the Green's function which satisfies (\ref{G}), see Section
\ref{Sect:2} below. 
For more distant $x$, $\bar \nu_{A,x}$ doesn't change
a lot any more: 
For $d\geq 2$, there exist constants $C_1(d)$ and $C_2(d)$
such that for all
$A\subset Q^d(n)$, $y\in A$, $x\in (U^d(an))^c$ with $a\geq 2\sqrt d$
\begin{equation} 
C_1 \bar \nu_{A,x} (y) \leq \bar \nu_{A,\infty} (y) \leq
C_2 \bar \nu_{A,x} (y), \label{eq2}
\end{equation}
see \cite{Lawlerbook}, Chapter 2.1. From (\ref{eq1}) and (\ref{eq2}), (A')
follows, and (A'') follows from (A') with (\ref{eq2}). 
Similarly we have the analogs of Theorem (B).
\endrem

(3) Our theorem improves a
result of Benjamini \cite{Benjamini}.
In fact, it implies the following weaker
statement (which is still stronger than 
\cite{Benjamini}): 
There exists $\rho (d)<d$, such that for any $\e >0$
there is an $N(\e)$  such that for any $n>N(\e)$, any $x \in
\Z^d$, and any $A\subset Q^d(n) = \{ 1,\dots , n\}^d$ there
is a set $\tilde A\subset A$ with 
$$ 
\omega (A, \tilde A, x) > \omega(A,A,x)-\e \quad \text{and} \quad
|\tilde A| <\e n^\rho .
$$ 
The analogous statements hold for harmonic measure conditioned on the
event that $A$ is reached, and also for harmonic measure from infinity.
Note that it is in
general impossible that
$\tilde A$ carries the full mass: 
Considering for example (for even $n$)
$A= \{ 1,3,5,\dots , n-1\}^d$, the only set having full mass (for
$x\not\in A$) is $A$, and
$|A|= (n/2)^d$. 
\endrem

(4) The dependence of the exponent $\rho$ on $\b$ for 2-dimensional
simple random walk paths $A$ (the ``multifractal spectrum
of the harmonic measure for $A$'') has been studied by 
Lawler \cite{Lawlerpreprint}.
Also for $d=2$, there
is another result of Lawler \cite{Lawler} which gives more information on
the support of harmonic measure from infinity $\bar\nu_{A,\infty}$ for
connected sets. 
\endrem

\section{Proof}

\subsection{Proof of Theorem (B)}

Take $n_K = 2^K$. Delete from $\{ 1,2,\dots ,n_K\}$ the 
central $\d 2^K$ points, from the remaining two
intervals of length $(1-\d) 2^{K-1}$ the central 
$\d (1-\d) 2^{K-1}$ points, and so on, $k$ $(<K)$ times.
In the $j$-th step, we have deleted $\d (1-\d)^{j-1} 2^{K-j+1}$ 
points and obtained intervals of length 
$(1-\d)^{j} 2^{K-j}$.
Let now $A_K$ be the product of $d$ copies of the resulting set. 
It consists of $2^{kd}$ squares of side length 
$(1-\d)^{k} 2^{K-k}$.
The total number of boundary points is
$$
| \partial A_K | = 2^{kd} \cdot 2d \cdot 
\left[ (1-\d)^{k} 2^{K-k}\right] ^{d-1}.
$$

To estimate the harmonic measure of the points of 
$\partial A_K$ we use the discrete Harnack inequality,
see for example
\cite{Lawlerbook}, Thm.\ 1.7.2: 
There exists a $c<\infty$ such that if 
$f: \Z^d \to [0,\infty)$ is
harmonic on $B_n$,
\begin{equation}\label{Harnack}
f(x_1)\leq c f(x_2), \quad \| x_1\|, \|x_2\| \leq n/2,
\end{equation}
with $B_n = \{ z\in\Z ^d\colon \| z\| <n\}$. 

Consider an arbitrary point $y\in \partial A_K$, 
and let $x_K$ be (for example) the central point of
$Q^d(n_K)$, i.e., $x_K= (2^{K-1}, \dots ,2^{K-1})$.
$Q^d(n_K)\setminus A_K$ consists of cylinders, 
called $j$-cylinders, of width $\d (1-\d)^{j-1} 2^{K-j+1}$,
$j=1,\dots , k$, in one component, and of width $n_K$ 
in the other components. 
$y$ lies on the boundary of a
$j_0$-cylinder for some $j_0\leq k$. 
Let $z_0$ be the point closest to $y$ lying in the 
center of the
$j_0$-cylinder. Let $z_1$ be the point closest to 
$z_0$ lying in the center of a $(j_0-1)$-cylinder. 
The distance from $z_0$ to $z_1$ is 
$\leq (1-\d)^{j_0-2}2^{K-j_0+1}$. 
Continue inductively to define points $z_i$
lying closest to $z_{i-1}$ in the center of a 
$(j_0-i)$-cylinder up to $i=j_0-1$. $|z_i - z_{i-1}|\leq
(1-\d)^{j_0-i-1}2^{K-j_0+i}$ and $|x_K - z_{j_0-1}|
\leq 2^{K-1}$. 
Applying (\ref{Harnack}) gives
\begin{eqnarray*}
\nu_{A_K,x_K} (y)& \geq &c^{-1/\d}\nu_{A_K,z_{j_0-1}} (y)
\geq c^{-1/\d}c^{-2/(\d(1-\d))}\nu_{A_K,z_{j_0-2}}
(y)\\
& \geq & \dots \geq c^{-1/\d} \left[c^{-2/(\d
(1-\d))}\right]^{j_0-1}\nu_{A_K,z_{0}}(y) \geq
c^{-4k/\d}\nu_{A_K,z_{0}}(y).
\end{eqnarray*}
We may estimate $\nu_{A_K,z_{0}}(y)$ simply by 
$\nu_{A_K,z_{0}}(y)\geq \tilde c \|z_0 - y\|^{1-d}\geq \tilde c
2^{-K(d-1)}$ (see \cite{Lawlerbook}, Lemma 1.7.4). 
Therefore
$$
\nu_{A_K,x_K} (y) \geq c^{-4k/\d}\tilde c
2^{-K(d-1)}.
$$

Now we want $| \partial A_K |>2^{K\rho}$ and 
$\nu_{A_K,x_K} (y) >2^{-K\b}$. 
This is achieved for large enough
$K$ by putting
$\d$ such that 
$\rho = d + 3(d-1)\log (1-\d) / \log 2$, 
$\b$ such that $\b-d+1$ $=$ $4 \log c/(\d \log 2)$, 
and $k=\g K$ with
$\g = \log \left[ 2(1-\d)^{3(d-1)}\right] $ 
$/\log \left[ 2(1-\d)^{d-1}\right]$.

\subsection{Discrete Hausdorff measure}

For bounded sets $A\subset \Z^d$, consider coverings of $A$ by a
countable number of balls $B_\a$ in $\Z^d$ with center $z_\a$ and radius
$r_\a$, $A\subset \Verein_\a B_\a$ with
$$ 
B_\a = \{ x\in \Z^d\colon \|x-z_\a \| \leq r_\a \} .
$$ 
For $0<\rho
\leq d$ we define
$$ 
h_\rho (A) = \inf \left\{ \sum_\alpha |B_\alpha |^{\rho / d} ; B_\alpha
\text{ ball}, A\subset \Verein_\a B_\a \right\} .
$$
Furthermore, consider a net of $l$-adic cubes: $\cC_0= \Z^d$, $\cC_1 =
\{ $cubes $C\subset \Z^d$ with side length $|C|^{1/d} = l$ and lower
corner $c=(k_1l, k_2l, \dots , k_dl)$ with $k_i\in \Z \}$,
$$ 
\cC_j =
\{ C\subset \Z^d\colon C=\{z\in\Z^d\colon k_il^j\leq z_i<( k_i+1)l^j
,k_i\in\Z, i=1\dots d\}\},
$$
and $\cC = \Verein_{j\in\N} \cC_j$.  
Analogously to $h_\rho$ we define
$$ 
m_\rho (A) = \inf \left\{ \sum_\alpha |C_\alpha |^{\rho / d} ; C_\alpha
\in\cC , A\subset \Verein_\a C_\a \right\} .
$$
Clearly, there exist two positive constants
$t_1(d)$ and $t_2(d,l,\rho)$ such that for all
$A\subset \Z^d$
\begin{equation}\label{hrhomrho}
 h_\rho (A) \leq t_1 (d) m_\rho (A)
\end{equation}
and
\begin{equation}\label{mrhohrho}
 m_\rho (A)\leq t_2(d,l, \rho) h_\rho (A) .
\end{equation}
By considering for example a ball of radius $\sqrt l$, one sees that the
dependence of
$t_2$ on $l$ cannot be removed. 
A possible choice is 
\begin{equation}\label{t2}
t_2= 8^d l^{d-\rho}.
\end{equation}

Analogously to Theorem 1  in Carleson \cite{Carleson}, p.7, (see also
\cite{Tsuji}, Chapter III.4) we have the following Lemma:
\begin{Lemma}\label{L:Carleson}
There  are constants $t_3$ and $t_4$, depending only on $d$, such that for
every bounded set $A\subset \Z^d$ there is a discrete measure
$\mu$ supported on $A$ with 
\begin{equation}\label{ballest}
\mu(B) \leq t_3|B|^{\rho /d} \quad \text{for all balls }B\subset \Z^d 
\end{equation}
and
\begin{equation}\label{massest}
\mu(A) \geq t_4 \, h_\rho(A).
\end{equation} 
\end{Lemma}

\noindent
{\bf Proof.} Start the construction of $\mu$ by putting $\mu_0(\{ x\} )=1$
for all $x\in A$ and $\mu_0(\{ x\})=0$ for
 $x\in A^c$. 
Choose your favorite $l$ and consider the cubes of
$\cC_1$. 
If for some
$C\in
\cC_1$
$\mu_0(C)>|C|^{\rho/d}$, reduce the density on the points of $C$
uniformly such that $\mu_1(C)=|C|^{\rho/d}$. 
Continue in this way. 
After
finitely many steps no further reduction will occur, since $\mu_k(C)\leq
|A|$ for all $C$ and $k$ and $|A| < l^{K\rho}$ for $K$ large enough. 
Put
$\mu=\mu_K$. 

$\mu$ satisfies 
$$
\mu (C)\leq |C|^{\rho/d} \quad \text{for all } C\in \cC 
$$
and therefore we have (\ref{ballest}).

From the construction of $\mu$, each point $a\in A$
is contained in a cube $C_\a$ with $\mu(C_\a) = |C_\a|^{\rho/d}$. 
If there
are several such cubes, choose the largest one. 
With this (disjoint)
covering $\{ C_\a\}$ we obtain
$$ 
\mu (A) = \sum_\a \mu (C_\a) = \sum_\a |C_\a|^{\rho/d}\geq m_\rho(A)
\geq \frac 1{t_1(d)}\, h_\rho(A)
$$ 
with (\ref{hrhomrho}).
This proves (\ref{massest}).\qed  

$\mu$ puts more mass on boundary points than on interior points. 
Thus it
is useful for estimating the harmonic measure, which is concentrated on
the boundary.

\subsection{Estimate of the trapping probability}\label{Sect:2}

Another useful quantity to estimate the harmonic measure in $d\geq 3$ is
the Green's function $G$, $G(x)$ being the expected number of visits to
$x$ of the random walk starting at 0,
$$
G(x) = \E ^0\left(\sum_{j=0}^\infty \charfct_{\{x\}}(S_j)\right) =
\sum_{j=0}^\infty\P^0 (S_j=x). 
$$
$G$ is harmonic in $\Z^d\setminus \{ 0\}$, $\Delta G(x)= -\delta (x)$,
and $G$ has the following asymptotic behavior:
$$ 
\lim_{\norm{x}\to \infty} \frac {G(x)}{a_d \norm{x}^{2-d}} =1, 
$$
where $a_d= 2/((d-2)\omega _d) $, and $\omega_d$ is the volume of the
unit ball in $\R^d$ (see for example \cite{Lawlerbook}, p.31). 
This
implies that there are constants $c_1$ and $c_2$ ($0<c_2<c_1$) depending
only on dimension such that we have the following upper and lower bounds
\begin{equation}  G(x) \leq c_1 \norm{x}^{2-d} \quad \text{and} \quad G(x)
\geq c_2 \norm{x}^{2-d} \quad \text{for}\quad x\in\Z^d\setminus \{ 0\} .
\label{G} \end{equation}

In $d=2$, $G$ is infinite, but there exists a quantity with similar
properties, namely the potential kernel
$$ 
a(x) = \lim_{n\to\infty} \sum_{j=0}^n\left( \P^0 (S_j=0) - \P^0
(S_j=x)\right).
$$
$\Delta a(x) = \delta(x)$, and $a$ has the following asymptotic behavior:
$$ 
\lim_{\norm{x}\to \infty} \left( a(x) - \frac 2\pi \log \norm{x}
-k\right) =0,
$$
where $k$ is some constant (see for example \cite{Lawlerbook}, p.38). 
Therefore there exists a constant $c$ such that we have the following upper and
lower bounds for $x\in\Z^d\setminus \{ 0\}$ 
\begin{equation} \label{a}  
a(x) \leq \frac 2\pi \log \norm{x}
+k +c \quad \text{and} \quad a(x)
\geq \frac 2\pi \log \norm{x}
+k - c .
\end{equation}

Consider now a cube $Q\subset \Z^d$, and let $Q_*\subset \Z^d$ be a cube
 of size $|Q_*|^{1/d}\leq q|Q|^{1/d}$,
where $q$ is a constant ($0<q<1$) to be determined below. 
$Q_*$ is placed
 such that its center is as close as possible to the center of $Q$: 
If
$|Q_*|^{1/d}$ and $|Q|^{1/d}$ are both even or both odd, $Q$ and $Q_*$ 
have the same center, and in the other cases, the distance of the centers
is $\sqrt d/2$. 
The next lemma gives an estimate of the probability that
a random walk starting in $Q_*$ reaches a set $A\subset \Z^d$ before
leaving $Q$, $\P ^a(\tau_A < \tau_{Q^c}) = \omega (A\verein Q^c,
A\schnitt Q, a)$:
\begin{Lemma}\label{L:B1}
Let $\rho \geq d-1$. 
Then for $q$ small enough (depending only on $d$) there exists $\tilde
c(d,q)>0$ such that for all $a\in Q_*$
\begin{equation}
\omega (A\verein Q^c, A\schnitt Q, a) \geq \tilde c \,
\frac {h_\rho(A\schnitt Q_*)}{|Q_*|^{\rho/d}} 
\label{omegaest}\end{equation} 
\end{Lemma}

\noindent
{\bf Proof.} If $A\schnitt Q_* = \emptyset$, (\ref{omegaest}) holds
trivially. 

Let now
$A\schnitt Q_* \neq \emptyset$ and let $\mu$ be the measure on $A\schnitt
Q_*$ from Lemma \ref{L:Carleson}. 
We treat first the case $d\geq 3$.
Consider the function
$u\colon \Z^d \to \R^+$,
$$ 
u (x) = \sum_{y\in A\schnitt Q_*} G(x-y)\, \mu(\{ y\}).
$$
$u$ is harmonic in $(A\schnitt Q_*)^c$. 
For $x\in Q_*$ and $y\in Q_*$,
$\norm{x-y}\leq |Q_*|^{1/d}\sqrt d$, and therefore with (\ref{G})
\begin{equation}\label{u*} 
u(x)\geq  c_2 d ^{(2-d)/2} |Q_*|^{(2-d)/d} \mu(A\schnitt Q_*) 
\quad  \text{for } x\in Q_* .
\end{equation}
For $x\in  Q^c$
and $y\in Q_*$,
$$
\norm{x-y}\geq  \frac{|Q|^{1/d} - |Q_*|^{1/d}}2  \geq
 \frac {1-q}{2q} |Q_*|^{1/d} 
$$ 
and therefore with (\ref{G})
\begin{equation}\label{ud} 
u(x)\leq  c_1 \left( \frac {1-q}{2q}\right)^{2-d}
|Q_*|^{(2-d)/d} \mu(A\schnitt Q_*)
\quad \text{for } x\in  Q^c .
\end{equation}
Furthermore, for all $x\in \Z^d$ 
\begin{equation}\label{us} 
u(x)\leq  c_3 |Q_*|^{(2+\rho -d)/d} ,
\end{equation}
where $c_3$ depends only on $d$.
This is seen as follows: First of all, with (\ref{G}),
$$
\sup_{x\in\Z^d} u(x) = \sup_{x\in  B(Q_*)} u(x) ,
$$
where $B(Q_*)$ is a ball with the same center as $Q_*$ and radius $a/2
\sqrt d |Q_*|^{1/d}$ with suitably chosen $a$ ($a =
1+2(c_1/c_2)^{1/(d-2)}$).
Now, for $x\in B(Q_*)$, 
$$ 
u(x) = \sum_{k=1}^{a \sqrt d |Q_*|^{1/d}} \sum _{y\in \tilde B_k(x)
}G(x-y)\,
\mu(\{ y\}) ,  
$$
where $\tilde B_k(x) = \{ y\in \Z^d\colon k-1\leq \| x-y\| <k\}$. 
Thus 
$$  
u(x) \leq G(0) \mu(\tilde B_1(x)) + \sum_{k=2}^{a \sqrt d |Q_*|^{1/d}}
c_1 (k-1)^{2-d} \mu (\tilde B_k(x)). 
$$
With $B_k(x)= \{ y\in \Z^d\colon  \| x-y\| <k\}$ we obtain
\begin{eqnarray*} 
 & & \sum_{k=2}^{a \sqrt d |Q_*|^{1/d}}
 (k-1)^{2-d} \mu (\tilde B_k(x)) \\
& & \qquad \qquad = (a \sqrt d |Q_*|^{1/d})^{2-d}
\mu(B_{a
\sqrt d |Q_*|^{1/d}}(x)) - \mu(B_1(x))\\
& &\qquad \qquad \quad + \sum_{k=2}^{a \sqrt d |Q_*|^{1/d}}
 \left( (k-1)^{2-d} - k^{2-d}\right) \mu (B_k(x)) .
\end{eqnarray*}
From (\ref{ballest}) we have $\mu(B_k(x)) \leq \tilde t_3 k^\rho$ for a
suitable $\tilde t_3 $ depending only on $d$. 
Then 
\begin{eqnarray*} 
\sum_{k=2}^{a \sqrt d |Q_*|^{1/d}}
 \left( (k-1)^{2-d} - k^{2-d}\right) \mu (B_k(x)) \leq t_3' \sum_{k=2}^{a
\sqrt d |Q_*|^{1/d}}
 k^{1-d+\rho}  & &\\
\leq t_3'\int_0^{a \sqrt d
|Q_*|^{1/d}+1} x^{1-d+\rho}dx =\frac{t_3'}{2-d+\rho} (a \sqrt d
|Q_*|^{1/d}+1)^{2-d+\rho}& &\\
\leq t_3' (a \sqrt d
|Q_*|^{1/d}+1)^{2-d+\rho}& &
\end{eqnarray*}
for $\rho\geq d-1$, where $t_3'$ depends only on $d$.
Putting everything together, we obtain (\ref{us}).

Consider now 
$$
\bar u(x) = \frac 1{\sup _{y\in \Z^d} u(y)} \left(u(x) - \sup _{y\in
Q^c} u(y)\right) .
$$
$\bar u(x)\leq 1$ for all $x\in \Z^d$, and $\bar u (x) \leq 0$ for $x\in
Q^c$. 
Compare $\bar u(x)$ with $\omega (A\verein Q^c, A\schnitt Q, x
) $: 
Application of the maximum principle (see for example
\cite{Lawlerbook}, p.25)  to $\bar u - \omega$ on $A^c\schnitt Q$ yields
$\bar u\leq \omega$ there, and on $A\schnitt Q$ we have $\omega =1\geq
\bar u$. 
Therefore
$$ 
\omega (A\verein Q^c, A\schnitt Q, x)\geq \bar u(x) \quad \text{for
all } x\in Q.
$$

Together with (\ref{u*}), (\ref{ud}), (\ref{us}), and (\ref{massest}), we
obtain for $a\in Q_*$
\begin{eqnarray*}
& & \omega (A\verein Q^c, A\schnitt Q, a)\\
& & \qquad \qquad \geq \frac
{\mu(A\schnitt Q_*)}{c_3 | Q_*|^{(2+\rho -d)/d}} \left( c_2 d ^{(2-d)/2}
-  c_1
\left(
\frac {1-q}{2q}\right)^{2-d}
\right) |Q_*|^{(2-d)/d} \\
& & \qquad \qquad \geq \tilde c \, \frac {h_\rho(A\schnitt
Q_*)}{|Q_*|^{\rho/d}} 
\end{eqnarray*} 
if we choose $q$ so small that $c_2 d ^{(2-d)/2}
-  c_1
\left(
({1-q})/{2q}\right)^{2-d}$ is
positive. 
This proves Lemma \ref{L:B1} in the case $d\geq 3$. 

For $d=2$, the analogous construction using instead of the Green's function
$G$ the potential kernel $a$ with
the estimates (\ref{a}) does the job.
\qed

Choose now $q$ so that Lemma \ref{L:B1} holds.

\subsection{An alternative for the cubes of $\cC$}

The estimate of the trapping probability (\ref{omegaest}) leads
to an alternative for the cubes $C$ of $\cC$: 
Either we have a local
estimate of the Hausdorff measure of $A\schnitt C$ or the harmonic
measure is localized on the outer shells of $C$. 
Cubes of the first kind
will be called (H)-cubes, those of the second kind (L)-cubes.

\begin{figure}
\begin{center}
\includegraphics[width=4cm]{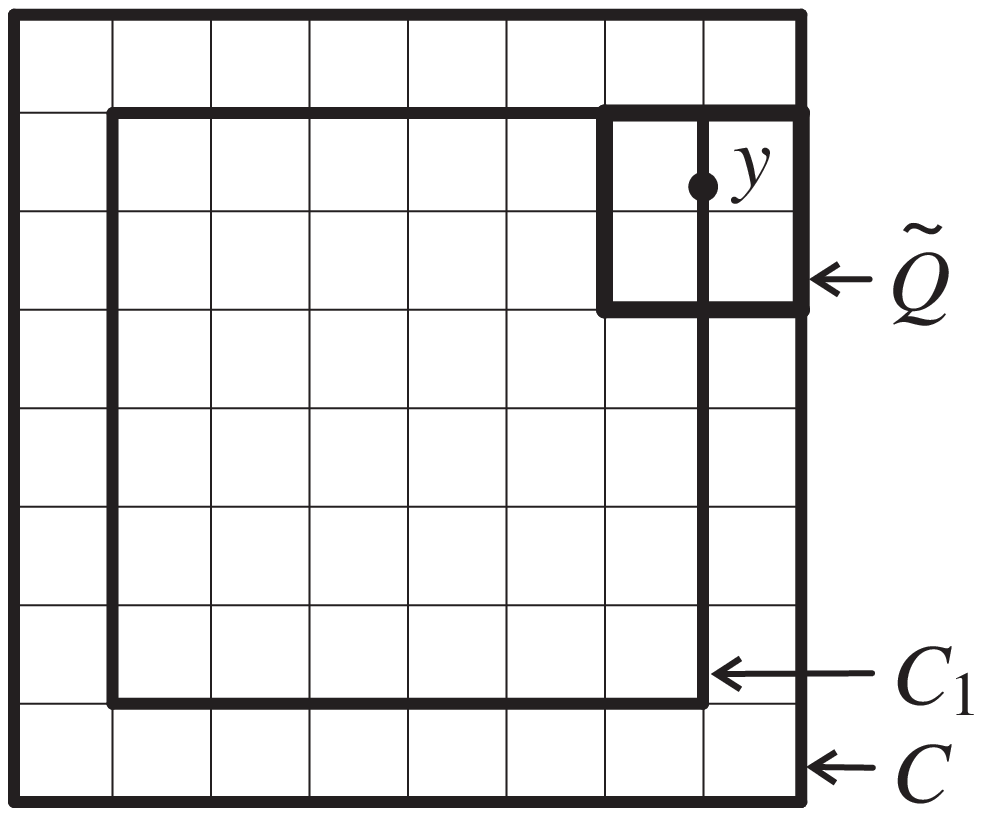}
\end{center}
\caption{For $d=2$ and $l=8$, this is a sketch of a cube $C\in \cC_j$
(for some $j$) together with its subcubes of $\cC_{j-1}$. 
By removing the
outer layer of subcubes, one obtains $C_1$. 
For $y\in \partial C_1$,
$\tilde Q$ is the union of the 4 nearest subcubes.}
\label{fig1}
\end{figure}

Consider now some $A\subset Q^d(n)$ and some $x\in \Z^d$. 
We abbreviate
$\omega(B)= \omega(A,B,x)$. 
For $C\in \cC_j$, $x \in (A\verein C)^c$, 
define (see
Fig.\ \ref{fig1})
\begin{eqnarray*}
C_1 & = &  C\setminus \text{outer subcubes }Q\in \cC_{j-1}, Q\subset C\\
C_2 & = &  C_1\setminus \text{outer }Q\text{'s in }C_1\\
\dots \quad C_{\bar l} & = & C_{\bar l-1}\setminus \text{outer }Q\text{'s
in }C_{\bar l-1}
\end{eqnarray*}
with $\bar l= l/6$. 
For  $x\in C\setminus A$, define
the $C_k$ by  successively removing layers of $Q$-cubes around the cube
$Q$ with $x\in Q$, and, if the boundary of $C$ is reached, remove also
successively layers of outer cubes like above. 

\begin{Lemma}\label{L:B2}
Let $\delta >0$ be small enough. 
Then for all $l$ there exists $\rho <d$ such that each cube
$C\in\cC_j$,
$j\geq 2$, satisfies one of the following conditions: 
\begin{eqnarray*} 
& \text{\rm{(H)}} & m_\rho (A\schnitt C) < |C|^{\rho
/d}\\ & \text{\rm{(L)}} &  
\omega( C_{\bar l}) \leq \frac {(1-c_4\d)^{\bar l-1}}{c_4\d}\ \omega(C),
\end{eqnarray*}
where $c_4$ is some constant depending only on $d$, $0<c_4<1$.
\end{Lemma}

\noindent
{\bf Proof.} Let $Q\in \cC_{j-1}$ be a subcube of $C$, and let 
$Q_*$ be the cube
 of size $|Q_*|^{1/d}= \left[ q|Q|^{1/d}\right]$ in the middle of $Q$.
From Lemma \ref{L:B1}, 
one of the following alternatives holds:
\begin{eqnarray}
& \omega (A\verein Q^c, A\schnitt Q, a) \geq \delta \quad \text{for all
}a\in Q_* & 
\label{Alt1} \\
& h_\rho(A\schnitt Q_*) < \displaystyle \frac \delta{\tilde c}\, 
|Q_*|^{\rho/d} & \label{Alt2}
\end{eqnarray}
We shall show that if (\ref{Alt1}) holds for all subcubes $Q\subset C$,
i.e., if we have a lower
bound for the trapping probability, then (L) holds for $C$, because the
harmonic measure will be concentrated on the outer shells. 
On the other
hand, if there is one subcube $Q$ with (\ref{Alt2}), we can estimate
$m_\rho$ of $A\schnitt C$. 

First case: 
There is a subcube $Q\subset C$, $Q\in \cC_{j-1}$, satisfying
(\ref{Alt2}). 
Then with (\ref{mrhohrho}), 
$$ 
m_\rho(A\schnitt Q_*) < \frac { t_2(d,l,\rho)\, \delta }{\tilde c }\, 
|Q_*|^{\rho/d} ,
$$
and
\begin{eqnarray*}
m_\rho (A\schnitt C) & \leq & m_\rho(C\setminus Q) + m_\rho(Q\setminus
Q_*) + m_\rho(A\schnitt Q_*) \\
 & \leq & (l^d-1)l^{(j-1)\rho} + l^d(1-q/2)^dl^{(j-2)\rho} + \frac
{t_2 \, \delta}{\tilde c } q^\rho l^{(j-1)\rho} .
\end{eqnarray*}
Now (H) follows if
\begin{equation}\label{C1}
l^d-1 + l^{d-\rho}(1-q/2)^d +
\frac {t_2(d,l,\rho) \, \delta}{\tilde c } q^\rho 
 < l^{\rho}  .
\end{equation}
Plug in (\ref{t2}) and choose
$\d$ so small that (\ref{C1}) for $\rho =d$ is satisfied, i.e., such
that $ (1-q/2)^d + 
8^d \delta q^d /\tilde c <1$. 
Then for all $l$ there exists $\rho <d$ such that (\ref{C1}) still holds.
Note that for large $l$ and small $d-\rho$,
(\ref{C1}) leads to 
\begin{equation}\label{dminusrho}
d-\rho \ {\approx} \ \frac b{l^d\log l}
\end{equation}
with $b= 1-[(1-q/2)^d + 8^d \delta q^d /\tilde c]$. 
We shall later choose $l$ very large and increasing with $\beta$.
Thus our $d-\rho$ goes to 0 as $\beta \to \infty$.

Second case: All subcubes $Q\subset C$, $Q\in \cC_{j-1}$, satisfy
(\ref{Alt1}). 
Since the
probability of running into $A$ before leaving
$Q$ is everywhere high, it is hard for the random walk to enter much into
the cube before having run into $A$, i.e., the harmonic measure of the
cubes deep inside $C$ will be very small. 
From the strong Markov property (see for
example \cite{Lawlerbook}, Theorem 1.3.2) we obtain
\begin{eqnarray*}
& & \omega( A\verein  C_k, C_k, x)  = \P^x (\tau_{A\verein 
C_k}<\infty ,\
S_{\tau_{A\verein  C_k}} \in C_k)\\
& &  = \sum_{y\in \partial C_{k-1}} \P^y (\tau_{A\verein 
C_k}<\infty ,\
S_{\tau_{A\verein  C_k}} \in C_k) \P^x (\tau_{A\verein 
C_{k-1}}<\infty ,\
S_{\tau_{A\verein  C_{k-1}}} =y) \\
& & \leq \sup_{y\in \partial C_{k-1}} \omega( A\verein 
C_k, C_k, y)
\ \omega( A\verein  C_{k-1}, C_{k-1}, x)
\end{eqnarray*}
(Here, $\partial A = \{ x\in A\colon \exists y\in A^c$ with $\| x-y\| =
1\}$.)   
Iterating this estimate, we get
\begin{equation}\label{Cbarl}
\omega( C_{\bar l}) \leq \omega ( A\verein C_{\bar l}, C_{\bar l}, x)
\leq \omega ( A\verein C_{1}, C_{1}, x) \prod_{k=2}^{\bar l}\sup_{y\in
\partial C_{k-1}} \omega( A\verein  C_k, C_k, y).
\end{equation}
On the other hand, using $\tau_{A\verein C_1}\leq \tau_A$ and the strong
Markov property, 
\begin{eqnarray}
\omega(C) & \geq & \sum_{y\in \partial C_1}\P^x (\tau_{A}<\infty ,\
S_{\tau_{A}} \in A\schnitt C,S_{\tau_{A\verein C_1}} =y ) \nonumber\\
& = & \sum_{y\in \partial C_1}\P^y (\tau_{A}<\infty ,\
S_{\tau_{A}} \in A\schnitt C) \P^x (\tau_{A\verein C_1}<\infty ,\
S_{\tau_{A\verein C_1}} =y )\nonumber\\
& \geq & \inf_{y\in \partial C_1} \omega (A, A\schnitt C,y)\ \omega (
A\verein C_{1}, C_{1}, x) \label{omegaC}
\end{eqnarray}
We shall show below that there exists a constant $c_4(d,q)$ such that
\begin{equation}
\omega (A, A\schnitt C,y) \geq c_4\d  \quad \text{for all } y\in \partial
C_1 , \label{c41}\end{equation}
and for $k=2,\dots ,\bar l$ 
\begin{equation}
\omega( A\verein  C_k, C_k, y) \leq 1-c_4\d  \quad \text{for all } y\in 
\partial C_{k-1} .
\label{c42}
\end{equation}
These estimates, together with (\ref{Cbarl}) and (\ref{omegaC}), yield
(L).

It remains to prove (\ref{c41}) and (\ref{c42}): 
Let $y\in \partial C_1$.
Consider the cube $\tilde Q$ formed from $2^d$ subcubes $Q\in\cC_{j-1}$ of
$C$ ``around'' 
$y$, i.e., the side length of $\tilde Q$ is $2l^{j-1}$, and the distance
of $y$ from the center of $\tilde Q$ is $\leq l^{j-1}/2 +1$ (see Fig.\
\ref{fig1}). 
We have $\tilde Q\subset C$, $\tilde Q \schnitt C_2 =
\emptyset$. 
Enumerate the $Q$-cubes in $\tilde Q$:
$$
\tilde Q = \Verein_{k=1}^{2^d}Q^{(k)} ,
$$
and let 
$$
\tilde Q_* = \Verein_{k=1}^{2^d}Q^{(k)}_* .
$$
Then, using again the strong Markov property, 
\begin{eqnarray*}
\omega (A, A\schnitt C,y) & = & \P^y (\tau_A<\infty , S_{\tau_A} \in
A\schnitt C)\\ & \geq & \P^y( \tau_{\tilde Q_*}< \tau_{\tilde
Q^c} \leq \infty , \exists t\in [\tau_{\tilde Q_*}, \tau_{\tilde
Q^c}) \text{ with }
S_t\in A) \\
& = & \sum_{a\in \tilde Q_*} \P^a(\tau_A < \tau_{\tilde Q^c}) \ \P^y (
\tau_{\tilde Q_* \verein \tilde Q^c} <\infty , S_{\tau_{\tilde Q_*
\verein \tilde Q^c}} =a)\\
& \geq & \sum_{k=1}^{2^d} \sum_{a\in \tilde Q_*^{(k)}} \P^a(\tau_A <
\tau_{ Q^{(k)c}}) \ \P^y (
\tau_{\tilde Q_* \verein \tilde Q^c} <\infty , S_{\tau_{\tilde Q_*
\verein \tilde Q^c}} =a)\\
& \geq & \d \ \omega (\tilde Q_*
\verein \tilde Q^c, \tilde Q_*, y),
\end{eqnarray*}
where we have used that all subcubes $Q\subset C$, $Q\in \cC_{j-1}$,
satisfy (\ref{Alt1}). 

\begin{figure}
\begin{center}
\includegraphics[width=4.7cm]{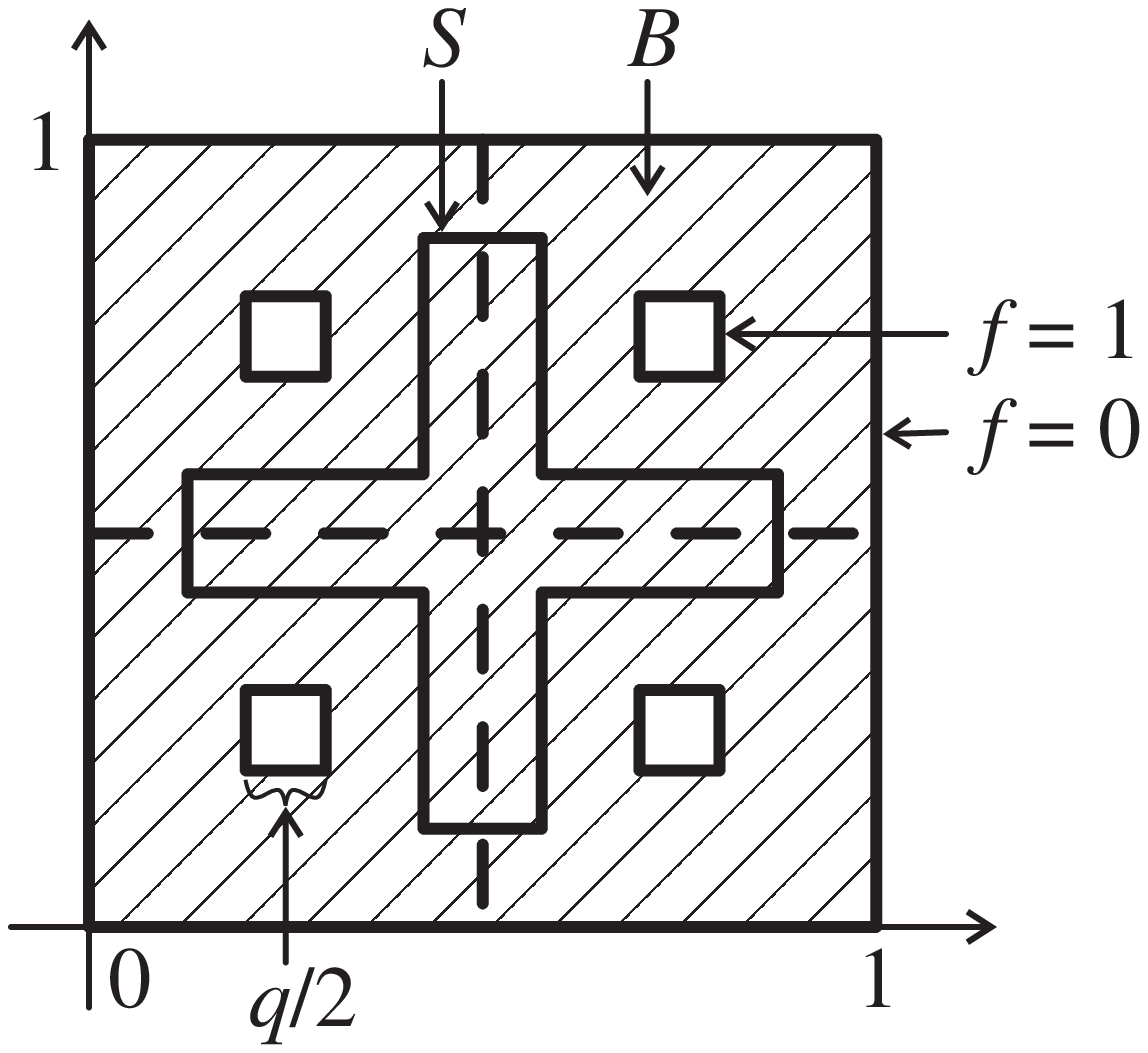}
\end{center}
\caption{For $d=2$, this is a sketch of the domain $B$ (hatched) $=
(0,1)^2 \setminus$ the 4 little squares of side length $q/2$. 
$B$
corresponds to $\tilde Q\setminus\tilde Q_*$, i.e., $(mB+z)\schnitt\Z^d$,
for suitable scale $m$ and shift $z$, equals $\tilde
Q\setminus\tilde Q_*$. 
The dashed middle axis lines correspond to the
boundaries of the subcubes making up
$\tilde Q$. 
The region $S$ is a neighborhood of the points
$x=m^{-1}(y-z)$ for those $y$'s which are possible for $\tilde
Q$, i.e., points on the middle half of a middle axis. 
The harmonic
function
$f$ on $B$ with boundary values
$f=0$ on the outer boundary of $B$ and $f=1$ on the boundaries of the
inner squares is bounded away from 0 on $S$.}
\label{fig2}
\end{figure}

To see that there exists $c_4$, independent of $l$
and $j$, with $\omega (\tilde Q_*
\verein \tilde Q^c, \tilde Q_*, y) \geq c_4$, remember that as a function
of $y$, $\omega (\tilde Q_*
\verein \tilde Q^c, \tilde Q_*, y)$ is (lattice) harmonic on $\tilde Q_*^c
\schnitt \tilde Q$ with boundary values $\omega = 1$ on $\tilde Q_*$ and
$\omega = 0$ on $\tilde Q^c$. 
Hence, the scaled function $\bar
\omega_m(x)= \omega (\tilde Q_*
\verein \tilde Q^c, \tilde Q_*, mx+z)$ with $m= 2l^{j-1}+1$, $2l^{j-1}$
the side length of $\tilde Q$, and suitable shift $z$, converges as
$m\to
\infty$ to the unique solution of $\Delta f=0$ on $B$, $f(x)=0$ on the
outer boundary of $B$ and $f(x)=1$ on the inner boundaries of $B$, where
$B$ is the ``limit'' of the scaled domains $m^{-1}(\tilde Q_*^c
\schnitt \tilde Q -z)$ as $m\to \infty$, see Fig.\ \ref{fig2}. 
Since the
convergence is uniform on compact subsets of $B$ \cite{convergence}, we
have a lower bound
$c_4$ for $\omega (\tilde Q_*
\verein \tilde Q^c, \tilde Q_*, y)$ for all $l$, $j$, and all $y=mx+z$
with $x$ in a region
$S$ around the middle halves of the middle axes of $B$ (see Fig.\
\ref{fig2}). 
This proves (\ref{c41}).

The proof of (\ref{c42}) is analogous: 
for $y\in \partial C_{k-1}$, put
$\tilde Q$ to be the cube consisting of $2^d$ subcubes of $C$ ``around''
$y$. 
Then $\tilde Q\schnitt C_k = \emptyset$ and $\tilde Q\subset C$. 
Thus
\begin{eqnarray*}
\omega (A\verein C_k,  C_k,y) & = & \P^y (\tau_{A\verein C_k}<\infty ,
S_{\tau_{A\verein C_k}}
\in  C_k)\\ 
& = &  \P^y (\tau_{A\verein C_k}<\infty ) - \P^y (\tau_{A\verein C_k}<\infty ,
S_{\tau_{A\verein C_k}}
\in A\setminus  C_k)\\
& \leq & 1- \P^y( \tau_{\tilde Q_*}< \tau_{\tilde
Q^c} \leq \infty , \exists t\in [\tau_{\tilde Q_*}, \tau_{\tilde
Q^c}) \text{ with }
S_t\in A) \\
& \leq & 1-c_4\d ,
\end{eqnarray*} 
with the same argument as above.
\qed

\subsection{Proof of Theorem (A)}\label{S:Proof}

Let now $\b>0$ and $n>N(\b)$ (to be chosen below). 
Let 
$A\subset Q^d(n)$, $x\in\Z^d$, and let $k^*\in \N$ be such
that $l^{k^*}\geq n>l^{k^*-1}$. 
To
the lower bound
$N(\b)$ there will correspond a $K^*$ such that $N(\b)=l^{K^*}$. 
We
construct Bourgain's tree $\cT$: starting with $C_0=\{ 1, \dots ,
l^{k^*}\}^d\in \cC_{k^*}$, we associate to each (L)-cube $C\in \cC_j$ its
$l^d$ subcubes in
$\cC_{j-1}$, and to each (H)-cube we associate a family $\{ C_\a\}$ with
$C_\a \subset C$, $A\schnitt C\subset \Verein _\a C_\a$, and $\sum_\a
|C_\a|^{\rho /d}<|C|^{\rho /d}$ (which exists according to Lemma
\ref{L:B2}). 
The elements of the tree are labeled by complexes $\g =
(\g_1, \dots , \g_k)$: $C_0$ has the label $\g =
(\g_1) = (0)$, its descendants have the label $\g =
(\g_1, \g_2) = (0, \g_2)$, and so on.

We stop the decomposition when  the cube is in
$\cC_1$ or
$\cC_0$ (because then Lemma
\ref{L:B2} doesn't apply any more). 
Thus each branch is at most $k^*$
long. 
Denote by $\g |k$ the restriction of $\g$ to the first $k$ digits. 
If $\tilde k$ is the length of $\g$, we call $C_{\g |1},C_{\g |2}, \dots,
C_{\g |\tilde k-1}$ the ``ancestors'' of $C_\g$. 
Let $\cT^*$ denote the set of the labels of the final cubes. 
We have
\begin{equation}\label{Aueberdeck}
A\subset \Verein_{\g\in \cT^*}C_\g .
\end{equation}

Given a maximal element $\g \in\cT^*$ of length $\tilde k$, we denote by
$\tau_k$ the length of the label of the $k$-th (L)-cube appearing in the
sequence
$C_{\g |1},C_{\g |2}, \dots$ of ancestors of $C_\g$, i.e., $C_{\g
|\tau_k}$ is the $k$-th (L)-cube, and 
$\tau_1<\tau_2<\dots < \tilde k$.
 ($\tau_k =\infty $
and $\g | \tau_k = \g$ if there are less than $k$ (L)-cubes in the
sequence $C_{\g |1},C_{\g |2}, \dots$ of ancestors of $C_\g$.) 

{\bf (a) Inner cubes.} The subcubes
$C_{\g
|\tau_k+1}$ of an (L)-cube $C_{\g
|\tau_k}$ are distinguished according to whether they lie in $(C_{\g
|\tau_k})_{\bar l}$ or not. 
If $x\in (A\verein C)^c$, the number of
subcubes which lie in $(C_{\g
|\tau_k})_{\bar l}$ is $(l-2\bar l)^d= (2/3)^d l^d$, and if 
$x\in C\setminus A$, the number  of
subcubes which lie in $(C_{\g
|\tau_k})_{\bar l}$ is simply estimated as $\geq (l-2\bar l)^d - (2\bar
l+1)^d \geq p l^d$ with $p= (2/3)^d - (1/2)^d$. 
To have a fixed proportion of
``inner'' subcubes (this simplifies the argument in part (c) below), we
shall choose for any (L)-cube $pl^d$ subcubes from those subcubes  
$C_{\g |\tau_k+1} \subset (C_{\g |\tau_k})_{\bar l}$ 
to call them ``inner'' subcubes. 

Let $ k_1 = k^*/3$ and $k_2=  (p/2) k_1$. 
Let 
\begin{eqnarray*}
\cT_< & = &  \{ \g\in \cT^* \colon \tau_{ k_1} (\g) =
\infty
\},\\
\cT_i & = & \{ \g\in \cT^* \colon \tau_{ k_1} (\g) < \infty ,\\ & & 
\text{ at least }k_2
\text{ of } C_{\g |\tau_1+1},C_{\g |\tau_2+1}, \dots,C_{\g
|\tau_{k_1}+1} \text{ are inner}\},
\end{eqnarray*}
and $\cT_{o} = \cT^* \setminus (\cT_< \verein \cT_i)$. 
If $C_{\g
|\tau_k+1}$ is inner, we have from Lemma \ref{L:B2}
$$ 
\omega (C_{\g
|\tau_k+1}) \leq  \omega \left( (C_{\g
|\tau_k})_{\bar l}\right) \leq \frac {(1-c_4\d)^{\bar l-1}}{c_4\d}\
\omega(C_{\g
|\tau_k}),
$$  
and if not, then in any case
$$ 
\omega (C_{\g
|\tau_k+1}) \leq  
\omega(C_{\g
|\tau_k}).
$$  
Then for $y\in \Verein_{\g\in \cT_i}C_\g$ we have (with $\g$ such that
$y\in C_\g$)
\begin{eqnarray*}
\nu_{A,x} (y)& \leq & \omega(C_\g) \leq \omega (C_{\g |
\tau_{k_1}+1})\leq 
\left( \frac {(1-c_4\d)^{\bar l-1}}{c_4\d}\right) ^{k_2}
\omega(C_{\g
|\tau_1})\\
& \leq &
\left( \frac {(1-c_4\d)^{\bar l-1}}{c_4\d}\right) ^{k_2}.
\end{eqnarray*}
Now choose $l$ such that
$$ 
\left( \frac {(1-c_4\d)^{\bar
l-1}}{c_4\d}\right) ^{k_2} < l^{-k^*\b } ,
$$
i.e.,
$$
 \frac p6 \left( \frac l6 -1\right)\log \frac 1{1-c_4\d} -\frac
p6
\log\frac 1{c_4\d} >
\b \log l.
$$
Then
$$
\Verein_{\g\in \cT_i}C_\g \subset \{ y\in \Z^d\colon \nu_{A,x} (y) <
n^{-\b}
\}.
$$
With (\ref{Aueberdeck}) we obtain
$$
\{ y\in\Z^d\colon \nu_{A,x} (y) \geq n^{-\b} \} \subset \Verein_{\g\in
\cT_<\verein \cT_o}C_\g .
$$
We shall show that $\sum_{\g \in \cT_<} |C_\g| \leq n^{-\tilde \rho}/2$
and $\sum_{\g \in \cT_o} |C_\g| \leq n^{-\tilde \rho}/2$ with 
$\tilde \rho = (\rho + d)/2$, where $\rho <d$ comes from Lemma \ref{L:B2}. 
This proves then Theorem (A).

{\bf (b) Estimate of $\sum_{\g \in \cT_<} |C_\g|$}. 
If $C_\g$ is
of type (H), then 
$$
\sum_{\g_k=1, \dots ; (\g, \g_k)\in \cT} |C_{(\g,
\g_k)}| ^{\rho /d} \leq |C_\g |^{\rho
/d},
$$
and if $C_\g \in \cC_j$ is of type (L),
then we have 
$$
\sum_{\g_k=1, \dots ,l^d} |C_{(\g, \g_k)}|
^{\rho /d} = l^{d}l^{(j-1)\rho }= l^{d-\rho} |C_\g |^{\rho
/d}.
$$
Thus 
\begin{eqnarray*}
|C_0|^{\rho/d} & \geq & \sum_{\g|\tau_1(\g);
\g\in\cT_<}|C_{\g|\tau_1(\g)}|^{\rho/d} \\
& \geq &l^{-(d-\rho)} \sum_{\g|\tau_1(\g)+1;
\g\in\cT_<}|C_{\g|\tau_1(\g)+1}|^{\rho/d} \\
& \geq & l^{-(d-\rho)} \sum_{\g|\tau_2(\g);
\g\in\cT_<}|C_{\g|\tau_2(\g)}|^{\rho/d} \\
& \geq & l^{-2(d-\rho)} \sum_{\g|\tau_2(\g)+1;
\g\in\cT_<}|C_{\g|\tau_2(\g)+1}|^{\rho/d} \\
\dots & \geq & l^{-(k_1-1)(d-\rho)}\sum_{\g\in \cT_<}
|C_\g|^{\rho /d}
\end{eqnarray*}
and therefore
\begin{equation*}
\sum_{\g \in \cT_<} |C_\g| = \sum_{\g \in \cT_<} |C_\g|^{\rho /d}
|C_\g|^{(d-\rho) /d} \leq l^{k_1(d-\rho)} l^{k^*\rho}.
\end{equation*}
For our choice of $ k_1$ and $\tilde \rho$ we have indeed
\begin{equation*} 
l^{k_1(d-\rho)} l^{k^*\rho} \leq
\frac 1 2\,  l^{(k^*-1)\tilde\rho}
\leq \frac1 2 \, n^{\tilde\rho}
\end{equation*}
for $k^*$ larger than some $K^*$.

{\bf (c) Estimate of $\sum_{\g \in \cT_o} |C_\g|$}.
Remember that $\cT_o = \{ \g\in \cT^* \colon \tau_{ k_1} (\g) < \infty ,$
$
\text{ less than }k_2
\text{ of } C_{\g |\tau_1+1},C_{\g |\tau_2+1}, \dots,C_{\g
|\tau_{k_1}+1} \text{ are inner}\}$. 
It is easy to see that 
\begin{equation*}
\sum_{\substack{\g\in \cT^*\colon \tau_{k_1}<\infty,\\
k \textrm{ of }  C_{\g |\tau_1+1},C_{\g |\tau_2+1}, \\ \dots,C_{\g
|\tau_{k_1}+1} \text{ are inner}}} |C_\g| \leq b(k;k_1,p)|C_0| = {k_1
\choose k} p^k (1-p)^{k_1-k}|C_0| ,
\end{equation*}
$b(k;k_1,p)$ being the binomial distribution, i.e., the distribution of
$\sum_{i=1}^{k_1}X_i$, where the $X_i$ are independent $\{ 0,1\}$-valued
random variables  with $P(X_i=1)=p$ for all $i$. 
For $0<a<p$, we have
from application of Markov's inequality to $\exp (\xi
\sum_{i=1}^{k_1}X_i)$
$$ 
P\left( \sum_{i=1}^{k_1}X_i \leq a k_1 \right)\leq e^{-k_1
I_p(a)}
$$ 
with 
$$
I_p(a) = a\log \frac ap + (1-a) \log \frac {1-a}{1-p}.
$$
(This is an elementary case of Cram\'er's theorem.) 
With $a= {k_2}/{k_1} =  p/2$,
$I_p(a)$ depends only on $d$. 
Then
\begin{equation*}
\sum_{\g \in \cT_o} |C_\g|  \leq  \sum_{j=0}^{k_2-1} b(j;k_1,p) |C_0|
\leq e^{-k_1 I} l^{k^*d},
\end{equation*}
and with our choice of the constants, noting also (\ref{dminusrho}), 
$$
e^{-k_1 I} l^{k^*d} \leq \frac 1 2\,  l^{(k^*-1)\tilde\rho}
\leq \frac1 2 \, n^{\tilde\rho} ,
$$
for $k^*$ larger than some $K^*$.
\qed

\end{document}